\documentclass[a4paper,12pt]{article}

\usepackage{setspace}
\usepackage{amsmath, amscd}
\usepackage{amssymb}
\usepackage{theorem}
\usepackage{color}
\usepackage{xypic}
\usepackage[english]{babel}
\usepackage[latin5]{inputenc}

\newtheorem{teo}{Theorem}[section]
\newtheorem{propo}[teo]{Proposition}
\newtheorem{lema}[teo]{Lemma}
\newtheorem{coro}[teo]{Corollary}

{\theorembodyfont{\rmfamily }
\newtheorem{defi}[teo]{Definition}

\newtheorem{ejem}[teo]{Example}}

\def\demo{\noindent \textit{Proof: }}

\def\K{\mathsf{K}}

\def\N{\mathtt{N}}
\def\T{\mathsf{T}}

\def\RR{\mathbb{R}}

\def\d{\mathrm{d}}

\def\dim{\mathrm{dim}\, }

\def\Hom{\mathrm{Hom}}

\def\Ker{\mathrm{Ker}\, }

\def\qed{\hspace*{\fill }$\square $}

\def\TenInv{\mathtt{T}}
\def\NatTen{\T }

\def\Sk{\mathsf{S}_{2\bar{p},k}}

\begin{document}

\title{Dimensional curvature identities \\ on pseudo-Riemannian geometry}

\date{November 15, 2013}

\author{Alberto Navarro \thanks{ICMat, Madrid, Spain.} 
 \and  Jos\'{e} Navarro \thanks{Department of Mathematics,
University of Extremadura, Avda. Elvas s/n, 06071, Badajoz, Spain. \newline {\it
Email address:} navarrogarmendia@unex.es \newline The second author has been partially supported by Junta de Extremadura and FEDER funds.}}

\maketitle

\begin{abstract}
For a fixed $n \in \mathbb{N}$, the curvature tensor of a pseudo-Riemannian metric, as well as its covariant derivatives, satisfy certain identities that hold on any manifold of dimension less or equal than $n$.

In this paper, we re-elaborate recent results by Gilkey-Park-Sekigawa  regarding these $p$-covariant  curvature identities, for $p=0,2$. To this end, we use the classical theory of natural operations, that allows us to simplify some arguments and to generalize the main results of Gilkey-Park-Sekigawa, both by dropping a symmetry hypothesis and by including $p$-covariant curvature identities, for any even $p$..

Thus, for any dimension $n$, our main result  describes the first space (i.e., that of highest weight) of $p$-covariant dimensional curvature identities, for any even $p$.



\end{abstract}

\tableofcontents



\section*{Introduction}

The curvature tensor of a pseudo-Riemannian metric, as well as its covariant derivatives, satisfy certain identities, such as the linear and differential Bianchi identities, or the Ricci identities. These identities are {\it universal}, in the sense that they are satisfied by the curvature tensor of {\it any} non-singular metric, on {\it any} manifold. Indeed, it can be proved that, essentially, these are the only identities with these properties (\cite{Kolar}).

Nevertheless, there exists some other kind of identities satisfied by the curvature.
As an example, recall that the Einstein tensor of a pseudo-Riemannian surface vanishes; that is, on any pseudo-Riemannian manifold of
dimension 2, the following
relation holds:
\begin{equation}\label{IdentidadEjemplo}
Ricc - \frac{r\, g}{2} \, = \, 0 \ ,
\end{equation} where $Ricc $ denotes the Ricci tensor and $r$ the scalar curvature of $g$.

This 2-covariant identity is satisfied by {\it any} non-singular metric on {\it any} manifold {\it of dimension less or equal} than two. For this reason, these type of identities will be called dimensional curvature identities (Definition \ref{DefinicionCentral}).


These dimensional identities already attracted attention very early in the development of General Relativity (\cite{Lanczos}). Later on, they also proved relevant in mathematics: in 1973, P. Gilkey  (\cite{Gilkey}) characterized the vanishing of the Pfaffian as the only scalar dimensional curvature identity satisfying certain homogeneity condition, and used this result to simplify the heat equation
proof of the index theorem. Let us also remark that, in the course of his investigations of the equivariant inverse problem of the calculus of variations, I. Anderson (\cite{Anderson}) was led to prove
a similar statement for symmetric, 2-covariant curvature identities, extending, in a sense, the identity (\ref{IdentidadEjemplo}) written above.

Nevertheless, most of these results were established in the realm of Riemannian manifolds, and their proofs relied on lengthy calculations, involving large expressions of multi-indexes. Recently, these arguments have been dramatically simplified by Gilkey-Park-Sekigawa (\cite{GilkeyI}), who also have extended their validity to pseudo-Riemannian manifolds of any signature (\cite{GilkeyII}), and have started studying analogous results in different settings (\cite{GilkeyIII}).

In this paper, we generalize the statements of Gilkey-Park-Sekigawa regarding scalar and symmetric, 2-covariant identities to the case of $p$-covariant identities, for any even $p$ and with no symmetry assumptions (Theorem \ref{Principal}). To this end, we use the machinery of the classical theory of natural operations (\cite{Kolar}, \cite{Tesis}), that also allows us to simplify the exposition of certain arguments.

The paper is organized as follows: the first Section is devoted to present the statement of the main result, Theorem \ref{Principal}, whereas the proofs are all postponed until the second Section.

\newpage

\section{Statements}

Fix a dimension $n \in \mathbb{N}$ and a signature $(n_+ , n_{-})$, where $n_+ + n_{-} = n$.

\begin{defi}\label{definicionNatTensor}
A  {\bf $p$-covariant, natural  tensor in dimension $n$}, {\bf associated to metrics of signature $(n_+ , n_{-})$}, is an assignment $T$ that, for any metric $g$ on any smooth manifold $X$ of dimension $n$, produces a $p$-covariant tensor $T(g)$ on $X$, satisfying:

\begin{itemize}
\item It is {\it regular}: 
if $\{ g_t \}_{t \in T}$ is a smooth family of metrics, parametrized by a smooth manifold $T$, then $\{ T(g_t) \}_{t\in T}$ is also a smooth family of tensors.

\item It is {\it natural}: for any local diffeomorphism $ \tau \colon Y \to X$ between smooth manifolds of dimension $n$, it holds
$$
T(\tau^* g)= \tau^* T (g) \ . $$
\end{itemize}

A natural tensor $ T $ is \textbf{homogeneous of weight} $w \in \mathbb{R}$ if, for any metric $g$ on any smooth manifold of dimension $n$, and any positive real number $\lambda>0$, it holds:
$$
T ( \lambda^2 g ) = \lambda^{w}\, T ( g)  \ . $$
\end{defi}

As it will be explained later (see Corollary \ref{NoImportaSignatura}), natural
tensors do not depend on the fixed signature of the metrics: the vector space of
homogeneous natural tensors associated to Riemannian metrics in dimension $n$ is
canonically isomorphic to the space of homogeneous natural tensors associated to
metrics of signature $(n_{+} , n_{-})$, provided that $n_{+} + n_{-} = n$.

Therefore, the $\RR$-vector space of homogeneous natural tensors will be denoted:
$$
\NatTen _{p,w} [n] :=   \left[
\begin{array}{c}
\text{ $p$-Covariant, natural tensors $T$ in dimension $n$ } \\
\text{ homogeneous of weight } w \
\end{array}
\right] \ . $$

\begin{ejem}\label{ExampleOne}
For any $n$, the metric itself $g$ is a natural 2-tensor in dimension $n$, homogeneous of weight 2.
The Riemann-Christoffel tensor, $R$, the Ricci tensor $Ricc$ and the scalar curvature $r$ are also natural tensors, homogeneous of weight 2, 0 and -2, respectively.
\end{ejem}

In general, it can be proved that the local expression of a homogeneous natural tensor is a ``universal'' (i.e., valid on any chart) polynomial function on the curvature and its covariant derivatives, with coefficients smooth functions of the metric and its inverse (\cite{Kolar}, \cite{Tesis}).

\subsection{Dimensional reduction. Universal tensors}

Let $(X,g)$ be an $(n-1)$-Riemannian manifold, an consider the cylinder $(X\times \RR \, , \, g + \d t^2)$, which is an $n$-dimensional Riemannian manifold.

Let $i$ denote the embedding:
$$ i \colon X \hookrightarrow  X\times \RR \qquad , \qquad x \ \mapsto \ (x,0) \ . $$

\begin{defi}
The \textbf{dimensional
reduction} of natural tensors is the linear map:
$$
\NatTen _{p,w} [n] \xrightarrow{\ \ r_n \ \ } \NatTen _{p,w} [n-1] \quad , \quad  r_n(T)(g):=i^* \left( T(g + \d t^2) \right) \ .
$$ \end{defi}

If $T \in \NatTen _{p,w} [n] $, it is not difficult to check that $r_n(T)$ is a natural tensor in dimension $n-1$. Moreover, as $\lambda^2 g + \d t^2$ and $\lambda^2 g + \lambda^2 \d t^2$ are related by an isometry that is the identity on $X$, $r_n (T)$ is also homogeneous of weight $w$:
\begin{align*}
r_n(T) (\lambda^2 g) &= i^* \left( T(\lambda^2 g + \d t^2 ) \right) = i^* \left( T ( \lambda^2 g + \lambda^2 \d t^2 ) \right) \\
& = i^* \left( \lambda^w T( g + \d t^2) \right) = \lambda^w \, r_n(T) ( g)  \ .
\end{align*}

Therefore, these linear maps establish a projective system:
$$
\ldots \xrightarrow{\quad r_{n+1} \quad} \NatTen _{p,w}[n] \xrightarrow{\quad r_n \quad}
\NatTen _{p,w}[n-1]  \xrightarrow{\quad r_{n-1}\quad} \ldots \ ,
$$
and the linear maps $r_n$ can be proved to be surjective (Lemma \ref{LemaUno}).

\begin{defi}
A {\bf universal tensor}, homogeneous of weight $w$, is an element of the inverse limit $$ \NatTen _{p,w}:= \varprojlim \NatTen _{p,w}[n] \ . $$
\end{defi}

\medskip

In other words, a universal tensor is a collection of natural
tensors $\, \{ T_n \}_{n\in \mathbb{N}} $,
where each $\,
T_{n}\,$ is a natural tensor in dimension $n$, satisfying that, for
every manifold $\,X\,$ of dimension $n_1$ and every embedding into a cylinder
$i \colon X \rightarrow X \times \mathbb{R}^{n_2}$, $x \ \mapsto \ (x,0)$, where $\mathbb{R}^{n_2}$ is endowed with the euclidean metric, it holds
$$i^*\left( T_{n_1 + n_2} ( g + \sum^{n_2}_{i=1} \d t_i^2 ) \right) = T_{n_1}(g) \ . $$

\medskip

\begin{ejem}
The metric $g$, the Riemann-Christoffel tensor $R$, the Ricci tensor $Ricc$, and the scalar curvature $r$ are all universal tensors.

For any fixed $\lambda\in \RR$, the tensor $\lambda g$ is  a universal tensor. However,  $ ({\rm tr} \, {\rm Id}) \, g$, $(-1)^{\dim X} g$ or $(-1)^{n_+} g$  are not universal tensors.


\end{ejem}

\medskip

Let $(X,g)$ be a Riemannian manifold as above and let $\pi \colon X \times \RR \to X$ be the first projection.

The curvature $R$ is a universal tensor that satisfies:
\begin{equation}\label{StrongProperty}
R( g + \d t^2) \, = \, \pi^* \left( R(g) \right) \ ,
\end{equation} which is a stronger property than that of being universal.

The following Lemma can be checked in local coordinates, using property
(\ref{StrongProperty}) above:

\begin{lema}\label{ContraccionCurvatura}
Let $T$ be a universal tensor, and let $T' = c ( R \otimes T)$ be the contraction of one index of the curvature with one index of $T$. Then $T'$ is a universal tensor.
\end{lema}

\begin{ejem}
Consider the metric $g$ as a 1-form with values on 1-forms, so that
$$ g^{2k + \bar{p}} := g \wedge \stackrel{2k + \bar{p}}{\ldots} \wedge g $$ is a $2k+ \bar{p}$-form with values on $2k + \bar{p}$-forms, for any $k, \bar{p}\geq 0$.

Analogously, consider the curvature tensor $R$ as a 2-vector with values on 2-vectors, so that:
$$ R^k := R \wedge \stackrel{k}{\ldots} \wedge R $$ is a $2k$-vector with values on $2k$-vectors.

Let $c$ be the contraction operator:
$$ \left( \Lambda^{2k} TX \otimes \Lambda^{2k + \bar{p} } T^* X \right) \otimes \left(
\Lambda^{2k} TX \otimes \Lambda^{2k + \bar{p} } T^* X \right)
\ \xrightarrow{ \quad c \quad }  \
\Lambda^{\bar{p} } T^* X \otimes \Lambda^{\bar{p} } T^* X \ . $$

For any $k, \bar{p} \geq 0$, let us define the $(2\bar{p})$-covariant universal tensors:
$$ \Sk := c \left( R^k \otimes g^{2k + \bar{p}} \right) \ . $$


Each $\Sk$ is indeed a universal tensor because it is obtained contracting indices of the curvature with indices of a universal tensor (Lemma \ref{ContraccionCurvatura}). Their local expression may be written as follows:
$$
(\Sk) _{i_1 \ldots i_{2\bar{p}}} :=
R^{a_1a_2,b_1b_2} \ldots R^{a_{2k-1}a_{2k}, b_{2k-1}b_{2k}} \delta_{b_1 \ldots
b_{2k} i_1 \ldots i_{\bar{p}}}^{c_1 \ldots c_{2k} j_1 \ldots j_{\bar{p}}} g_{a_1 c_1} \ldots g_{a_{2k} c_{2k}} g_{j_1 i_{\bar{p} +1}} \ldots g_{j_{\bar{p}} i_{2\bar{p}}}
$$
where $\delta_{i_1 \dots i_m}^{j_1 \ldots j_m}$ denotes the generalized Kronecker delta.

Moreover, these tensors $\Sk$ are skew-symmetric in its first $\bar{p}$
indexes and in the last $\bar{p}$ indexes, and, due to the symmetries of $g$ and $R$, they are also symmetric under the interchange of these group of indexes.

If $\bar{p}=0$, there is only one $\Sk$, whereas, for $\bar{p} \geq 1$, any permutation $\sigma$ of $2\bar{p}$ elements produces the tensor
$$ (\sigma \cdot \Sk ) (D_1, \ldots , D_{2\bar{p}}) \, := \, \Sk (D_{\sigma (1)} , \ldots ,
D_{\sigma (2\bar{p})} ) \ . $$

Due to the symmetries of $\Sk$, there may only be $\frac{1}{2} \left( p \cdot
(p-1) \cdot \ldots \cdot (\frac{p}{2}+1) \right)$ different tensors among the $\{
\sigma \cdot \Sk \}_{\sigma \in S_{2\bar{p}} }$.
\end{ejem}

\subsection{Dimensional curvature identities}

Loosely speaking, any natural tensor is locally written in terms of
the coefficients of the curvature and its covariant derivatives. Hence, if, for some $n$,  a universal tensor lies in the kernel of the canonical projection $\NatTen _{p,w}  \longrightarrow \NatTen _{p,w} [n] $, then it can be understood as an identity, satisfied by the coefficients of the curvature and its covariant derivatives, which is valid for any metric, of any signature, on any manifold of dimension less or equal than $n$.

Of course, if a universal tensor $\{ T_n \}_{n \in \mathbb{N}}$ defines a curvature identity, then the universal tensor $\{
\lambda T_n \}_{n \in \mathbb{N}}$ defines the same identity, for any $\lambda \in \RR -\{ 0 \}$.

This motivates the following definition:

\begin{defi}\label{DefinicionCentral}
A \textbf{dimensional curvature identity} in dimension $n$ is an element of the projective space associated to the vector space:
$$ \K _{p , w} [n] := \Ker \left[ \phantom{\frac{1}{1}} \hskip -.3cm  \NatTen _{p,w}  \longrightarrow \NatTen _{p,w} [n] \right]   \ . $$
\end{defi}

\begin{ejem}\label{SkSonIdentidades} For any fixed $\bar{p}, k \geq 0$ (apart from the cases $k = 0$ and $\bar{p} = 0,1$) the tensor $\Sk$ vanishes whenever $\dim X < 2k + \bar{p}$, because the form $g^{2k + \bar{p}}$ is identically zero.

Hence, this tensor $\Sk$ defines a dimensional curvature identity
$$  \Sk \, \in \, \K_{2 \bar{p}, w} \left[ 2k + \bar{p} - 1 \right] \ , $$
with $w = 2 (\bar{p} - k)$.

The following Theorem \ref{Principal} establishes that it is, essentially, the only dimensional curvature identity of this kind.
\end{ejem}

Finally, let us observe that, if $T \in \NatTen _{2 \bar{p} , w} $ is a
$(2\bar{p})$-covariant, universal tensor, homogeneous of degree $w$, then $w$ has to
be an even integer, lesser or equal than $p$ (see, v. gr., \cite{EinsteinNS}), so
that we may write, without loss of generality:
$$ w = 2 (\bar{p} - k) $$ for some $k \geq 0$.

\begin{teo}\label{Principal}

Consider covariant tensors with an even number $2\bar{p}$ of indices ($\bar{p} \geq 0$), associated to non-singular metrics.

For any\footnote{Apart from the exceptional cases $(p,k)=(0 ,0)$ or $(1,0)$} weight $w = 2 (
\bar{p} - k) $, with $k \geq 0$, the following holds:

\begin{itemize}
\item If $n \geq 2k + \bar{p}$, there are no dimensional curvature identities of weight $w$; that is,
$$  \K_{2\bar{p} , w} \left[ n \right] \, = \, 0 \quad , \quad  \mbox{ for } n \geq 2k + \bar{p} \ . $$

\item The vector space $\K_{2\bar{p} , w} \left[ 2k + \bar{p} - 1 \right] $ is generated by the tensors $\{ \sigma \cdot \mathsf{S}_{k} \}_{\sigma \in S_{2\bar{p}}}$, and hence has dimension:
$$ \dim \left( \K_{2\bar{p} , w} \left[ 2k + \bar{p} - 1 \right] \right) \, =  \, \bar{p} \, \left( \bar{p} + 1 \right) \cdot \ldots \cdot (2\bar{p}-2)  (2\bar{p}-1)   \  , $$
or reduces to a single identity, when $\bar{p}=0$.
\end{itemize}
\end{teo}

In the particular case $\bar{p}=0$ and $n = 2k$, the smooth function $\mathsf{S}_{0 ,k}$ is proportional to the Pfaffian of the curvature.  Hence, Theorem \ref{Principal} implies the following result:

\begin{coro}[\cite{GilkeyII}]
Consider tensors with $p=0$ indices; i.e., scalar differential invariants associated to non-singular metrics.

For any weight $w = -2k$, with $k \geq 1$, the following holds:

\begin{itemize}
\item If $n \geq 2k$, there are no dimensional curvature identities of weight $w$; i.e.:
$$ \K_{0 ,w} \left[ n \right] \, = \, 0 \quad , \quad \mbox{ for } n \geq 2k  \ .  $$

\item If $n = 2k-1$, the only dimensional curvature identity of weight $w = -2k$ is the vanishing of $\mathsf{S}_{0 ,k} $:
$$ \K _{0 ,w} \left[ 2k-1 \right] \, = \, \langle \mathsf{S}_{0 ,k} \rangle \ . $$
\end{itemize}
\end{coro}


In the particular case $\bar{p}=1$, the symmetric 2-tensors producing the identities are well-known: $\mathsf{S}_{2,1}$ is the Einstein tensor; $\mathsf{S}_{2,2}$ was first introduced by Lanczos (\cite{Lanczos}) and the other $\mathsf{S}_{2, k} $ were independently introduced by Lovelock (\cite{LovelockI}, \cite{LovelockAl}) and Kuz'mina (\cite{Kuzmina}). 

Theorem \ref{Principal} then reads:

\begin{coro}
Consider covariant tensors with $p=2$ indices (not necessarily symmetric) associated to non-singular metrics.

For any weight $w = 2 -2k$, with $k \geq 1$, the following holds:

\begin{itemize}
\item If $n \geq 2k + 1$, there are no dimensional curvature identities of weight $w$; i.e.:
$$ \K_{2,w} \left[ n \right] \, = \, 0 \quad , \quad \mbox{ for } n \geq 2k + 1 \ .  $$

\item If $n = 2k $, the only dimensional curvature identity of weight $w= 2 - 2k$ is the vanishing of $\mathsf{S}_{2, k} $:
$$ \K_{2,w} \left[ 2k \right] \, = \, \langle\mathsf{S}_{2, k}  \rangle \ . $$
\end{itemize}
\end{coro}

This statement drops the symmetry hypothesis that is assumed in \cite{GilkeyII}, as well as in the Riemanian statements of \cite{Anderson} and \cite{GilkeyI}.

\section{Proofs}\label{proofs}

The proof of Theorem \ref{Principal} relies on the classical theory of natural
constructions, that reduces the problem to a question regarding tensors (at a point)
invariant under the action of the orthogonal group (see Theorem
\ref{StredderSlovak}).

\medskip

To state this result, let us firstly introduce the vector spaces of normal tensors:

\begin{defi}\label{definitionnormal} Let $X$ be a smooth manifold of
dimension $n$, $x \in X$ be a point and $r\geq 2$ be an integer. The space  $N_r
\subset S^2 T^*_xX \otimes S^r T^*_xX $ of  {\bf$r^{th}$-order metric normal tensors}
at  $x$ is the kernel of the symmetrization in the last $(r+1)$-indices:
$$ 0 \to N_r \to S^2 T^*_xX \otimes S^r T^*_xX \xrightarrow{\ s_{r+1} \ } T^*_xX \otimes S^{r+1} T^*_xX \to 0 \ . $$
\end{defi}

\medskip
Any germ of metric $g$ around the point $x$ defines a sequence of metric normal tensors $(g^2_x , \ldots, g^r_x , \ldots ) \in N_2 \times \ldots \times N_r \times \ldots$; to construct this sequence, choose normal coordinates $x_1 , \ldots , x_n$ for $g$ at $x$ and define
$$ (g^k_x)_{ab i_1 \ldots i_k} \, := \, \frac{\partial^k g_{ab}}{\partial x_{i_1} \ldots \partial x_{i_k}} (x) \quad , \qquad k = 2 , \ldots  $$ for the condition of the chart $(x_i)$ being normal guarantees that the symmetrization of the last $r+1$ indices of $g^r_x$ is zero.

Let $\NatTen _{p,w} [n_+ , n_{-} ]$ denote the vector space of homogeneous natural
$p$-tensors (of weight $w$), associated to pseudo-Riemannian metrics of signature
$(n_+ , n_{-})$.

\begin{teo}[\cite{Tesis}, \cite{Stredder}]\label{StredderSlovak}
Let $X$ be a smooth manifold of dimension $n$, $x\in X$ be a point and $g_x$ be a pseudo-Riemannian metric at $x$ of signature $(n_+ , n_{-})$.

There exists
an $\mathbb{R}$-linear isomorphism:
$$
\NatTen _{p,w} [n_+ , n_{-} ] \ \simeq \  \bigoplus \limits_{d\in D} \, \mathrm{Hom}_{O_{g_x}} \left(
S^{d_2}N_2 \otimes \cdots \otimes S^{d_r}N_r   \ , \ \otimes^p T_x^* X \
\right)
$$
where $D$ is the set of sequences of nonnegative integers $d=\{d_2 , \ldots , d_r \}
$ such that:
\begin{equation}\label{Condicion}
2d_2 + \ldots + r\, d_r = p - w \ .
\end{equation}
If such equation has no solutions the such vector space is zero.
\end{teo}

\medskip

If $\varphi \colon S^{d_2} N_2 \otimes \cdots \otimes S^{d_r} N_r \to \otimes^p T_x^*
X$ is an $O_{g_x}$-equivariant linear map, then, on any metric $g$ with the prefixed
value at $x$, the corresponding natural tensor $T$ is obtained by the formula:
$$ T(g)_{x} = \varphi \left( (g^2_{x} \otimes \stackrel{d_2}{\ldots} \otimes \, g^2_{x}) \otimes \cdots
\otimes (g^s_{x} \stackrel{d_s}{\ldots} \otimes \, g^s_{x}) \right) $$ where $(g^2_{x}, g^3_{x}, \ldots )$ is the sequence of normal tensors of $g$ at the point $x \in X$.

The value at any other point, and for any other metric over another manifold, is computed adequately transforming with a diffeomorphism.

\medskip

Indeed, the $O_{g_x}$-equivariant linear maps in the theorem can, in certain cases, be explicitly computed
applying the invariant theory of the orthogonal group explained in the next section.
An interesting consequence of this Theorem \ref{StredderSlovak} and Corollary \ref{NoImportaSign} below is that the vector spaces $\NatTen _{p,w} [n]$ do not depend on the signature of the metrics under consideration:

\begin{coro}\label{NoImportaSignatura}
If $n_+ + n_{-} = m_+ + m_{-}$, there is a canonical isomorphism:
$$\, \NatTen _{p,w} [n_+ , n_{-} ] = \NatTen _{p,w} [m_+ ,m_{-} ]\, . $$
\end{coro}



\subsection{Invariant theory for the orthogonal group}

Let $(E,g)$ be a $\mathbb{R}$-vector space of dimension $n$ with a nonsingular metric
$g$ of signature $(n_+,n_{-})$ and let $O_g$ denote the Lie group of its linear isometries
$(E,g) \to (E,g)$.

The main theorem of the invariant theory for the orthogonal group describes the
polynomial functions on $m$ vectors
$$
f \colon E \times \stackrel{m}{\ldots} \times E \longrightarrow \mathbb{R}
$$ that are invariant under the action of $O_g$.

For any given $i,j = 1, \ldots , m $, the following functions $y_{ij}$ are examples of
$O_g$-invariant polynomial functions:
$$
y_{ij} \colon E \times \stackrel{m}{\ldots} \times E \to \mathbb{R} \quad , \quad
y_{ij} (e_1, \ldots , e_m) := g (e_i , e_j) \ .
$$


If $m > n := n_{+} + n_{-}$,  then such functions have relations: for any $1 \leq i_0 < \ldots < i_n
\leq m$, $1 \leq j_0 < \ldots < j_n \leq m$, the following identities hold:
\begin{equation*} 
\left| \begin{matrix}
y_{i_0j_0} & \ldots & y_{i_0 j_n} \\
\vdots &  & \vdots \\
y_{i_nj_0} & \ldots & y_{i_nj_n}
\end{matrix}
\right| \, (e_1, \ldots , e_{m})\,   = (e_{i_0} \wedge \ldots \wedge e_{i_n} ) \cdot
( e_{j_0} \wedge \ldots \wedge e_{j_n}) = 0 \cdot 0 = 0  ,
\end{equation*}
where $\cdot$ denotes the metric induced by $g$ on the corresponding tensor algebra.

The so-called Main Theorem then states that these are, essentially, the only invariant functions and the only relations among them:


\begin{teo} [\cite{Carlos}, \cite{Weyl}]\label{Main}

The algebra $\mathsf{A}_m^g$ of $O_g$-invariant polynomial functions on $E \times \stackrel{m}{\ldots}
\times E $ is generated by the functions $y_{ij}$.

Moreover, let $Y_{ij}$ be free symmetric variables. The map $Y_{ij}\mapsto y_{ij}$ induces a canonical isomorphism
$$\begin{CD}
  \mathbb{R} [Y_{ij}] / M_{n+1} @=  \mathsf{A}_m^g  \ ,
\end{CD}$$
where $M_{n+1}$ is the ideal generated by the functions:
\begin{equation*}\label{ThisSetting}
 M_{i_0 \ldots i_n}^{j_0 \ldots j_n} = \left| \begin{matrix}
Y_{i_0j_0} & \ldots & Y_{i_0 j_n} \\
\vdots &  & \vdots \\
Y_{i_nj_0} & \ldots & Y_{i_nj_n}
\end{matrix} \right| \ , \ \mbox{ for any } \ \ \begin{matrix}
1 \leq i_0 < \ldots < i_n \leq m\\
1\leq j_0 < \ldots < j_n \leq m
\end{matrix} \ .
\end{equation*}


In particular, if $m \leq n$, these functions $y_{ij}$ are algebraically independent.
\end{teo}

\medskip

As a consequence, if $(E' , g')$ is another $n$-dimensional $\RR$-vector space with a nonsingular metric of signature $(n'_{+} , n'_{-})$, then there are canonical isomorphisms:
$$
\begin{CD}
 \mathsf{A}^{g'}_{m} @= \mathbb{R} [Y_{ij}] / M_{n+1} @=  \mathsf{A}_m^g  \ ,
\end{CD}$$

This classical statement is usually proved in the realm of algebraic varieties (\cite{Carlos}, \cite{Weyl}); that is, for the {\it affine algebraic $\mathbb{R}$-group} $O_g$. The corresponding version for the (non-compact) Lie group $O_g$ requires some argument, to reduce the proof to the algebraic case (see, for example, \cite{Jaime}, \cite{GilkeyII}, or \cite{Tesis}).

As a consequence, it readily follows a useful description of the space of linear
forms
$$
E \otimes \stackrel{m}{\ldots} \otimes E \xrightarrow{ \quad  \quad }
\mathbb{R}
$$
which are invariant under the action of $O_g$:

\begin{coro}\label{TensoresInvariantes}
The vector space $\Hom _{O_g} (E \otimes \buildrel{m}\over{\ldots} \otimes E,
\mathbb{R}) $ of invariant linear forms is zero if $m$ is odd and, if $m = 2k$ is
even, it is spanned by total contractions:
\begin{equation*}\label{Generadores}
\omega_\sigma \colon e_1
\otimes \ldots \otimes e_{2k} \ \mapsto \ g(e_{\sigma(1)} , e_{\sigma (2)}) \cdot
\ldots \cdot g(e_{\sigma (2k-1)} , e_{\sigma (2k)})
\end{equation*}
where $\sigma\in S_{2k}$ is a permutation.

Moreover, if $m \leq 2n$, the only relations among these generators are the obvious ones due to the symmetry of $g$.
\end{coro}

\demo Among the polynomials on $y_{ij}$, observe that $m$-multilinear maps are
precisely the linear combinations of
$$
y_{\sigma(1) \sigma(2)} \cdot  \ldots  \cdot y_{\sigma(2k-1) \sigma(2k)}
$$
where $\sigma$ is a permutation of $1, \ldots , 2k $.





\qed

\medskip

The isomorphisms $\mathsf{A}^{g'}_{m} = \mathbb{R} [Y_{ij}] / M_{n+1} =  \mathsf{A}_m^g $ takes $m$-multilinear maps into $m$-multilinear maps. Hence:

\begin{coro}\label{NoImportaSign} Let $(E,g)$, $(E', g')$ be vector spaces of the same dimension,
 endowed with non-singular metrics. There exists a canonical isomorphism:
$$
\Hom _{O_g} (E \otimes \buildrel{m}\over{\ldots} \otimes E, \mathbb{R}) \ = \ \Hom
_{O_{g'}} (E' \otimes \buildrel{m}\over{\ldots} \otimes E' , \mathbb{R}) \ .
$$
\end{coro}

From now on, we put
$$ \TenInv _m [n] := \mathrm{Hom}_{O_g} ( E \otimes \buildrel{m}\over{\ldots} \otimes E , \mathbb{R} ) \ $$ where $g$ is any non-singular metric; for example a scalar product.




\begin{defi}
For all $n > 1$, there exist \textbf{dimensional reduction} linear maps
$$
r_n \colon \TenInv _m [n] \longrightarrow \TenInv _m [n-1]  \ ,
$$ defined as follows: if $\omega \in\TenInv _m [n]$ is an invariant linear form on $n$-dimensional euclidean spaces, let $r_n(\omega)$ be the invariant form on $(n-1)$-dimensional euclidean spaces obtained as:
$$
 E \otimes \ldots \otimes E \hookrightarrow (E \perp \mathbb{R} )
\otimes \ldots \otimes ( E\perp \mathbb{R}) \xrightarrow{\ \omega \ } \mathbb{R} \ . $$
\end{defi}

The following statement is a consequence of  Theorem \ref{Main}:

\begin{propo}\label{LemaSimplifica}
The dimensional reduction maps $\, r_n \colon \TenInv _m [n] \longrightarrow \TenInv
_m [n-1]\,$ are surjective for all $n > 1$ and, for $n
> m -1$, they are linear isomorphisms.
\end{propo}

\subsubsection{Invariant forms on the subspace of normal tensors}

Let $X$ be a smooth manifold of dimension $n$, $g_x$ a Riemannian metric at a point $x \in X$, and let $D = (d_2, \ldots , d_r)$ be a multi-index.


As the orthogonal group is semisimple, restriction to the spaces of normal tensors induce surjective linear maps, for each $n \in \mathbb{N}$:
\begin{equation*}\label{Extension}
\xymatrix{
\TenInv _D[n] := \mathrm{Hom}_{O_{g_x}} \left(
T_x^* X \otimes \stackrel{4d_2 + \ldots + (2+r)d_r + p}{\ldots} \otimes T^*_xX \ , \ \mathbb{R}
\right) \, \ar[d]_-{i^*_n}  \\
\N_D [n] := \mathrm{Hom}_{O_{g_x}} \left(
S^{d_2}N_2 \otimes \cdots \otimes S^{d_r}N_r \otimes^p T_x^* X \ , \ \mathbb{R}
\right) \ . }  
\end{equation*} 



That is to say,

\begin{lema}\label{Semisimplicidad}
Any $O_{g_x}$-invariant linear map:
$$
S^{d_2}N_2 \otimes \cdots \otimes S^{d_r}N_r   \otimes^p T^*_x X \ \longrightarrow \
\mathbb{R} \
$$
is the restriction of a $O_{g_x}$-invariant linear map:
\begin{equation*}
T^*_x X \otimes \ \stackrel{4d_2 + \ldots + (2+r)d_r + p}{\ldots} \ \otimes T^*_xX
\ \longrightarrow \ \mathbb{R} \ .
\end{equation*}
\end{lema}





In terms of the generators $\omega_\sigma$ introduced in Corollary \ref{TensoresInvariantes}, the dimensional reduction maps $r_n$ have a simple expression:
$$ \TenInv _D[n] \, \xrightarrow{\quad r_n \quad } \TenInv _D[n-1] \qquad , \qquad \omega_\sigma \ \longmapsto \ \omega_\sigma \ . $$

Therefore, they specialize to the subspaces of normal tensors, defining maps $\bar{r}_n$:
$$
\xymatrix{
\TenInv _D[n] \, \ar[r]^-{r_n} \ar[d]_-{i^*_n }  & \TenInv _D[n-1] \ar[d]^{i^*_{n-1}}  \\
\N_D [n] \ar@{.>}[r]^-{\bar{r}_n} & \N_D [n-1] }
$$

The following Lemma states that any relation satisfied by the dimensional reduction $\bar{r}_n$ in the subspace 
is indeed the restriction of a relation satisfied by the dimensional reduction $r_n$ in the ambient space: 

\begin{lema}\label{LemaTecnico} Restriction to the space of normal tensors induce surjective linear maps:
$$
\xymatrix{
\mathtt{K}_D[n] := \Ker \left( r_n \colon \TenInv _D[n] \, \to \TenInv _D[n-1] \right) \, \ar[d]_-{i^*_n}  \\
\bar{\mathtt{K}}_D[n] := \Ker \left( \phantom{\frac{1}{1}} \hskip -.2cm \bar{r}_n \colon \N _D[n] \, \to \N _D [n-1] \right) } $$
\end{lema}

\demo The kernel of the restriction maps $i^*_n$
is generated by those total contractions $\omega_\sigma$ that become zero when restricted to $ S^{d_2} N_2 \otimes \ldots \otimes S^{d_r} N_r \otimes^p T^*_xX$. Due to the particular (``geometrical'') definition of these subspaces, the vanishing of these contractions does not depend on the dimension $n$, but only on the symmetries defining the normal tensors and the symmetric powers. Hence, the dimensional reductions $r_n$ induce a surjective maps
\begin{equation*}
 \Ker i^*_n \, \xrightarrow{\ r_n \ } \, \Ker i^*_{n-1}  \ .
\end{equation*}

Now, the thesis follows using the Snake's Lemma.

\qed

\subsection{Final computations}

A system of generators for $\T_{2\bar{p},w} \left[ n \right]$ is given, via Theorem \ref{StredderSlovak}, by total contraction maps
$$ \omega_\sigma \colon S^{d_2} N_2 \otimes \ldots \otimes S^{d_r} N_r \otimes^{2\bar{p}} T^*_xX \longrightarrow
\mathbb{R} $$ where $\sigma \in S_m$ is a permutation of $m = 4d_2 + \ldots + (2+r)d_r + 2\bar{p}$ indices, and the multi-index $D = (d_2, \ldots , d_r)$ satisfies:
$$ 2 d_2 + \ldots + r d_r \, = \, 2\bar{p} - w \ . $$


In terms of these generators, the dimensional reduction maps are:
$$ \NatTen _{2\bar{p},w}
[n] \xrightarrow{\ \ r_n \ \ } \NatTen _{2\bar{p},w} [n-1] \quad , \quad \omega_\sigma \, \longmapsto \, \omega_\sigma \ .  $$

\begin{lema}\label{LemaUno}
For any weight $w = 2\bar{p} - 2k$, with $k \geq 0$, the dimensional reduction maps
$$ r_n \colon \NatTen _{2\bar{p},w} [n] \xrightarrow{\ \  \ \ } \NatTen _{2\bar{p},w} [n-1]  $$
satisfy:

\begin{itemize}
\item They are surjective, for all $n$.

\item If $n > 2k + \bar{p}$, they are linear isomorphisms:

\item If $n = 2k + \bar{p}$, then any tensor on the kernel of $r_n$ is second-order.
\end{itemize}
\end{lema}

\demo Theorem \ref{StredderSlovak}, in conjunction with Proposition \ref{LemaSimplifica}, imply that $r_n$ is surjective, for all $n$.

To study the kernels, first observe that the number $m =  4d_2 + \ldots + (2+r)d_r + 2\bar{p}$ of indices to contract on a generator $\omega_\sigma$ is bounded by:
\begin{align*}
m &= \ 4 d_2 + \ldots + (2 + r) d_r + 2\bar{p}  \ = \ 2 ( d_2 +
\ldots + d_r ) +  2 d_2 + \ldots + rd_r  + 2\bar{p}   \\
&  \leq  \ 2\bar{p} - w + 2\bar{p} - w + 2\bar{p} \ =  \ 4k + 2\bar{p} \ = \ 2  \left( 2k + \bar{p} \right) \ .
\end{align*}

Therefore, if the dimension is big enough, $n > 2k + \bar{p}$, then there are no relations among the generators $\omega_\sigma$ due to dimensional considerations (Corollary \ref{TensoresInvariantes}).

Nevertheless, on a manifold of dimension $n = 2k + \bar{p} -1$, it may happen that the (even) number of indices $m$ is strictly greater than twice the dimension; in that case,
\begin{align*}
2 \left( 2k + \bar{p} \right) \ & \leq \  m = 4 d_2 + \ldots + (2 + r) d_r + 2\bar{p} \ \leq  \  2 \left( 2k + \bar{p} \right) \ ,
\end{align*}
so both inequalities are indeed equalities, and hence $d_3 = \ldots
= d_r = 0$.

This amounts to saying that the corresponding universal tensor is second-order.

\qed

\begin{lema}\label{LemaDos}
For any weight $w = 2\bar{p} -2k$, with $k \geq 0$, it holds:
$$ \dim \, \K _{2\bar{p} , w} \left[ 2k + \bar{p} -1  \, \right] \, = \,  \bar{p} \, \left( \bar{p} + 1 \right) \cdot \ldots \cdot (2\bar{p}-2)  (2\bar{p}-1) \ .  $$

If $p=0$, then $\dim \K _{0,w} [2k -1] = 1$, for all $k \geq 0$.
\end{lema}

\demo By the previous Lemma, any tensor in the kernel of $r_n$, for $n = 2k + \bar{p}$, is second-order,  so that it is defined by a linear combination of invariant linear maps:
$$
\omega \, \colon \, S^{\frac{n}{2}} N_2 \otimes^{2\bar{p}} T^*_x X
\longrightarrow \mathbb{R} \ ,
$$
which are non-zero if $\dim X =  n + \bar{p}$, but vanish when $\dim
X = n + \bar{p} - 1$.

Theorem \ref{Main}, together with Corollary \ref{TensoresInvariantes} and Lema \ref{LemaTecnico}, imply that any such $\omega$  is a linear combination of elements of the form
\begin{equation}\label{Generador}
e_1 \otimes \ldots \otimes e_{2(n+\bar{p})} \  \longmapsto \
(e_{r_{1}} \wedge \ldots \wedge e_{r_{n+\bar{p}}} ) \cdot ( e_{s_1} \wedge \ldots
\wedge e_{s_{n+\bar{p}}})
\end{equation} where the indexes $r_{1}, \ldots , r_{n+\bar{p}}, s_1 , \ldots , s_{n+\bar{p}}$ run from 1 to $2(n + \bar{p})$.

Let us now prove that, due to the symmetries of the space $S^{\frac{n}{2}} N_2 \otimes^{2\bar{p}} T^*_x X$, we can extract, among these generators, the following basis:
\begin{align*}\label{LaUnica}
\omega_{i_1 \ldots i_{\bar{p}}}^{j_1 \ldots j_{\bar{p}}} \ \colon \ S^{\frac{n}{2}} N_2 &\otimes^{2\bar{p}} T^*_x X  \longrightarrow \mathbb{R} \\
e_1 \otimes \ldots \otimes e_{2(n+\bar{p})}   \ \longmapsto \ (e_1 \wedge \ldots
 \wedge e_{2n-1}\wedge e_{i_1} & \wedge  \ldots \wedge e_{i_{\bar{p}}}  ) \cdot (e_2
\wedge \ldots \wedge e_{2n}\wedge e_{j_1} \wedge \ldots \wedge
e_{j_{\bar{p}}} )
\end{align*} where $i_1, \ldots , i_{\bar{p}} , j_1 , \ldots , j_{\bar{p}}$ are different indexes, running from $2n+1$ to $2(n+\bar{p})$.

So let $\omega$ be a non-zero linear map as in (\ref{Generador}). Up to a sign, we can assume  that $e_{r_1} = e_1 $. Since normal tensors in
$N_2$ are symmetric in the first two indexes, we may also assume that $e_{s_1} = e_2 $.
Analogously, normal tensors of order two are
symmetric on the third and four indexes, so that we can also write $e_{r_2}=e_3$ and
$e_{s_2}=e_4$.

A similar argument easily proves that $\omega$ is proportional to the linear map that sends $e_1 \otimes \ldots \otimes e_{2(n+\bar{p})}$ into
$$
(e_{1} \wedge e_3 \wedge \ldots \wedge e_{2n-1}\wedge e_{i_1}\wedge \ldots \wedge
e_{i_{\bar{p}}} ) \cdot (e_{2} \wedge e_4 \wedge \ldots \wedge e_{2n}\wedge
e_{j_1}\wedge \ldots \wedge e_{j_{\bar{p}}} )
$$ where $i_1, \ldots , i_{\bar{p}} , j_1 , \ldots , j_{\bar{p}}$ are different indexes, running from $2n+1$ to $2(n+\bar{p})$; that is to say, $\omega$ is proportional to the $\omega_{i_1 \ldots i_{\bar{p}}}^{j_1 \ldots j_{\bar{p}}} $ defined above.

If $\bar{p}=0$, we are done. Otherwise, the dimension of the kernel of $r_n$ is bounded by the number of possible unordered choices of $\bar{p}$ elements over a set $2\bar{p}$ elements, and divided by 2 (because the metric in $\Lambda^{k} E$ is symmetric). In other words, it is bounded by
$$
\bar{p} \, \left( \bar{p} + 1 \right) \cdot \ldots \cdot (2\bar{p}-2)  (2\bar{p}-1) \ .
$$

As the linear maps $\omega_{i_1 \ldots i_{\bar{p}}}^{j_1 \ldots j_{\bar{p}}} $  are linearly independent, the statement follows.

\qed

\medskip

\noindent {\it Proof of Theorem \ref{Principal}:} As it was explained in Example \ref{SkSonIdentidades}, the universal tensors $\sigma \cdot \Sk$ define elements in $\K_{2 \bar{p}, 2 (\bar{p} - k)} \left[ 2k + \bar{p} - 1 \right]$, for any permutation $\sigma$ of $2\bar{p}$ elements.

The tensors $\sigma \cdot \Sk$ are all $\mathbb{R}$-linearly independent, so Theorem \ref{Principal} readily follows from Lemma \ref{LemaUno} and Lemma \ref{LemaDos}

\hfill $\square$

\bigskip

\subsection*{Acknowledgements}

The authors thank  J. A. Navarro and J. B. Sancho for their generous advice and helpful comments.

The second author has been partially supported by Junta de Extremadura and
FEDER funds.

\end{document}